\documentclass[11pt,fullpage]{article}

\newcommand{\ol}{\setlength{\itemsep}{0pt.}\begin{enumerate}}
\newcommand{\eol}{\end{enumerate}\setlength{\itemsep}{-\parsep}}
\newcommand{\ignore}[1]{}
\title{\bf Hyperbolic Polynomials Approach to Van der Waerden/Schrijver-Valiant like Conjectures :\\
Sharper Bounds , Simpler Proofs and Algorithmic Applications}
 
\author{Leonid Gurvits \thanks{%
{\tt gurvits@lanl.gov}. Los Alamos National Laboratory, 
Los Alamos, NM. } 
}

\begin{document}

%\sloppy

\begin{titlepage}

\maketitle

\begin{abstract}
Let $p(x_1,...,x_n) = p(X) , X \in R^{n}$ be a homogeneous polynomial of degree $n$ in $n$ real variables ,
$e = (1,1,..,1) \in R^n$ be a vector of all ones . Such polynomial $p$ is
called $e$-hyperbolic if for all real vectors $X \in R^{n}$ the univariate polynomial
equation $p(te - X) = 0$ has all real roots $\lambda_{1}(X) \geq ... \geq \lambda_{n}(X)$ .
The number of nonzero roots $|\{i :\lambda_{i}(X) \neq 0 \}|$ is called $Rank_{p}(X)$ .
A $e$-hyperbolic polynomial $p$ is called $POS$-hyperbolic if roots
of vectors $X \in R^{n}_{+}$ with nonnegative coordinates are also nonnegative
(the orthant $R^{n}_{+}$ belongs to the hyperbolic cone) and $p(e) > 0$ .
Below $\{e_1,...,e_n\}$ stands for the canonical orthogonal basis in $R^{n}$. \\
The main results states that if $p(x_1,x_2,...,x_n)$ is a $POS$-hyperbolic (homogeneous) polynomial of degree $n$ ,
$Rank_{p} (e_{i}) = R_i$  and 
$
p(x_1,x_2,...,x_n) \geq \prod_{1 \leq i \leq n} x_i ; x_i > 0 , 1 \leq i \leq n ,
$ \\
then the following inequality holds
$$
\frac{\partial^n}{\partial x_1...\partial x_n} p(0,...,0) \geq  \prod_{1 \leq i \leq n}  (\frac{G_{i} -1}{G_{i}})^{G_{i} -1}  (G_i = \min(R_{i} , n+1-i)) .
$$

This theorem is a vast (and unifying) generalization of the van der Waerden conjecture on the permanents of doubly stochastic matrices as well as
the Schrijver-Valiant conjecture on the number of perfect matchings in $k$-regular bipartite graphs .
These two famous results correspond to the $POS$-hyperbolic polynomials being
 products of linear forms.   
 
 Our proof is relatively simple and "noncomputational" ; it actually slightly improves Schrijver's lower bound ,
 and uses very basic ( more or less centered around
 Rolle's theorem ) properties of hyperbolic polynomials . \\
We present some important algorithmic applications of the result, including a polynomial time deterministic algorithm approximating 
the permanent of $n \times n$ nonnegative entry-wise matrices within a multiplicative factor $\frac{e^n}{n^m}$ for any fixed positive $m$ . 
This paper introduces a new powerful "polynomial" technique , which allows as to simplify/unify famous and hard known results
as well to prove new important theorems .

The paper is (almost) entirely self-contained , most of the proofs can be found in the {\bf Appendices}.

\end{abstract} 
\end{titlepage}
\newpage

%  THEOREM-LIKE ENVIRONMENTS
 
\newtheorem{THEOREM}{Theorem}[section]
\newenvironment{theorem}{\begin{THEOREM} \hspace{-.85em} {\bf :} 
}%
                        {\end{THEOREM}}
\newtheorem{LEMMA}[THEOREM]{Lemma}
\newenvironment{lemma}{\begin{LEMMA} \hspace{-.85em} {\bf :} }%
                      {\end{LEMMA}}
\newtheorem{COROLLARY}[THEOREM]{Corollary}
\newenvironment{corollary}{\begin{COROLLARY} \hspace{-.85em} {\bf 
:} }%
                          {\end{COROLLARY}}
\newtheorem{PROPOSITION}[THEOREM]{Proposition}
\newenvironment{proposition}{\begin{PROPOSITION} \hspace{-.85em} 
{\bf :} }%
                            {\end{PROPOSITION}}
\newtheorem{DEFINITION}[THEOREM]{Definition}
\newenvironment{definition}{\begin{DEFINITION} \hspace{-.85em} {\bf 
:} \rm}%
                            {\end{DEFINITION}}
\newtheorem{EXAMPLE}[THEOREM]{Example}
\newenvironment{example}{\begin{EXAMPLE} \hspace{-.85em} {\bf :} 
\rm}%
                            {\end{EXAMPLE}}
\newtheorem{CONJECTURE}[THEOREM]{Conjecture}
\newenvironment{conjecture}{\begin{CONJECTURE} \hspace{-.85em} 
{\bf :} \rm}%
                            {\end{CONJECTURE}}
\newtheorem{PROBLEM}[THEOREM]{Problem}
\newenvironment{problem}{\begin{PROBLEM} \hspace{-.85em} {\bf :} 
\rm}%
                            {\end{PROBLEM}}
\newtheorem{QUESTION}[THEOREM]{Question}
\newenvironment{question}{\begin{QUESTION} \hspace{-.85em} {\bf :} 
\rm}%
                            {\end{QUESTION}}
\newtheorem{REMARK}[THEOREM]{Remark}
\newenvironment{remark}{\begin{REMARK} \hspace{-.85em} {\bf :} 
\rm}%
                            {\end{REMARK}}
\newtheorem{FACT}[THEOREM]{Fact}
\newenvironment{fact}{\begin{FACT} \hspace{-.85em} {\bf :} 
\rm}%                     
		            {\end{FACT}}

%\newenvironment{proof}{\noindent {\bf Proof:} \hspace{.677em}}%
%                      {}
 
%theorem
\newcommand{\thm}{\begin{theorem}}
%lemma
\newcommand{\lem}{\begin{lemma}}
%proposition
\newcommand{\pro}{\begin{proposition}}
%definition
\newcommand{\dfn}{\begin{definition}}
%remark
\newcommand{\rem}{\begin{remark}}
%example
\newcommand{\xam}{\begin{example}}
%conjecture
\newcommand{\cnj}{\begin{conjecture}}
%problem
\newcommand{\prb}{\begin{problem}}
%question
\newcommand{\que}{\begin{question}}
%corollary
\newcommand{\cor}{\begin{corollary}}
%fact
\newcommand{\fac}{\begin{fact}}

%proof
\newcommand{\prf}{\noindent{\bf Proof:} }
%end theorem
\newcommand{\ethm}{\end{theorem}}
%end lemma
\newcommand{\elem}{\end{lemma}}
%end proposition
\newcommand{\epro}{\end{proposition}}
%end definition
\newcommand{\edfn}{\bbox\end{definition}}
%end remark
\newcommand{\erem}{\bbox\end{remark}}
%end example
\newcommand{\exam}{\bbox\end{example}}
%end conjecture
\newcommand{\ecnj}{\bbox\end{conjecture}}
%end problem
\newcommand{\eprb}{\bbox\end{problem}}
%end question
\newcommand{\eque}{\bbox\end{question}}
%end corollary
\newcommand{\ecor}{\end{corollary}}
%end fact
\newcommand{\efac}{\end{fact}}
%end proof
\newcommand{\eprf}{\bbox}
%begin equation
\newcommand{\beqn}{\begin{equation}}
%end equation
\newcommand{\eeqn}{\end{equation}}
% white box
\newcommand{\wbox}{\mbox{$\sqcap$\llap{$\sqcup$}}}
%black box
\newcommand{\bbox}{\vrule height7pt width4pt depth1pt}
\newcommand{\qed}{\bbox}

%right arrow
\newcommand{\rarrow}{\rightarrow}
%left arrow
\newcommand{\larrow}{\leftarrow}
%gradient
\newcommand{\grad}{\bigtriangledown}

\overfullrule=0pt
\def\setof#1{\lbrace #1 \rbrace}
\section{Introduction}
An $n \times n$ matrix $A$ is called doubly stochastic if it is nonnegative entry-wise
and every column and row sum to one. The set of $n \times n$ doubly stochastic
matrices is denoted by $\Omega_{n}$. Let $\Lambda(k,n)$ denote the set of
$n \times n$ matrices with nonnegative integer entries and row and column sums all
equal to $k$ . We define the following subset of rational doubly stochastic matrices :
$\Omega_{k,n} = \{ k^{-1} A : A \in \Lambda(k,n) \}$ .
In a 1989 paper \cite{bapat}  R.B. Bapat defined the set 
$D_{n}$ of doubly stochastic $n$-tuples of $n \times n$ matrices.\\  
An $n$-tuple ${ \bf A }=(A_{1},
\cdots, A_{n})$ belongs to $D_{n}$ iff 
$A_{i} \succeq 0$, i.e. $A_{i}$ is a positive
semi-definite matrix, $1 \leq i \leq n$ ;
$tr A_{i} = 1$ for $1 \leq i \leq n$ ;
$\sum^{n}_{i=1} A_{i} = I$, where $I$, as usual,
stands for the identity matrix.  
Recall that the permanent of a square matrix A is defined by
$$ per(A) = \sum_{\sigma \in S_{n}} \prod^{n}_{i=1} A(i,
\sigma(i)). $$
Let us consider an $n$-tuple
${ \bf A } = (A_{1}, A_{2},... A_{n})$, where $A_{i} = (A_{i}(k,l):
1 \leq k, l \leq n)$ is a complex $n \times n$ matrix $(1 \leq i
\leq n)$.  Then $\det (\sum_{1 \leq i \leq n} t_{i}A_{i})$ is a homogeneous
polynomial of degree n in $t_{1},t_{2}, \cdots, t_{n}$.  The number 
\begin{equation}
M({ \bf A }):=D(A_{1}, A_{2}, \cdots, A_{n}) =
\frac{\partial^{n}}{\partial t_{1} \cdots \partial t_{n}} \det
  (t_{1} A_{1} + \cdots + t_{n}A_{n}) 
\end{equation}
is called the mixed discriminant of $A_{1}, A_{2}, \cdots,
A_{n}$.\\
The permanent is a particular (diagonal) case of the mixed discriminant . I.e.
define a multilinear polynomial $Mul_{A}(t_1,...,t_n) = \prod_{1 \leq i \leq n} \sum_{1 \leq j \leq n} A(i,j)t_j$.
Then $per(A) = \frac{\partial^{n}}{\partial t_{1} \cdots \partial t_{n}} Mul_{A}(t_1,...,t_n).$

Let us recall two famous results and one recent result by the author.
\begin{enumerate}
\item {\bf Van der Waerden Conjecture}\\
The famous Van der Waerden Conjecture \cite{minc} states that
$ min_{A \in \Omega_{n}} D(A) = \frac{n!}{n^{n}}$ {\bf (VDW-bound)} and the minimum is
attained uniquely at the matrix $J_{n}$ in which every entry equals
$\frac{1}{n}$. Van der Waerden Conjecture was posed in 1926
and proved only in 1981 : D.I. Falikman proved
in \cite{fal} the lower bound
$\frac{n!}{n^{n}}$ ; 
the full conjecture , i.e. the uniqueness part , was proved by G.P. Egorychev in \cite{ego} .

\item {\bf Schrijver-Valiant Conjecture}\\
Define 
$$
\lambda(k,n) = \min \{per(A) : A \in \Omega_{k,n}\} = k^{-n} \min \{per(A) : A \in \Lambda_{k,n}\} ; 
\theta(k) = \lim_{n \rightarrow \infty}(\lambda(k,n))^{\frac{1}{n}} .
$$
It was proved in \cite{schr-val} that , using our notations , $\theta(k) \leq g(k) = (\frac{k-1}{k})^{k-1}$ and conjectured
that $\theta(k) = g(k)$ . Though the case of $k=3$ was proved by M. Voorhoeve in 1979 \cite{vor} , this conjecture was settled
only in 1998 \cite{schr} (17 years after the published proof of the Van der Waerden Conjecture). The main result of \cite{schr}
is the following remarkable inequality :\\
$\min \{per(A) : A \in \Omega_{k,n}\} \geq (\frac{k-1}{k})^{(k-1)n}$ {\bf (Schrijver-bound)} . \\
The proof in \cite{schr} is probably
one of the most complicated and least understood in the theory of graphs .
\item {\bf Bapat's Conjecture (Van der Waerden Conjecture for mixed discriminants) }\\
One of the problems posed in \cite{bapat} is to determine the minimum
of mixed discriminants of doubly stochastic tuples :
$ min_{A \in D_{n}} D(A) = ?$ \\ 
Quite naturally, R.V.Bapat conjectured that
$ min_{A \in D_{n}} D(A) = \frac{n!}{n^{n}}$ {\bf (Bapat-bound)}
and that it is attained uniquely at ${\bf J}_{n} =:(\frac{1}{n} I,...,
\frac{1}{n} I)$.\\
In \cite{bapat} this conjecture was formulated for real matrices. The author
had proved it \cite{gur} for the complex case, i.e. when matrices
$A_{i}$ above are complex positive semidefinite and, thus,
hermitian.
\end{enumerate}
The {\bf (VDW-bound)}
is the simplest and most powerful bound on permanents and therefore among the simplest
and most powerful general purpose bounds in combinatorics. Besides its many applications to
the graph theory and combinatorics , {\bf (VDW-bound)} has been recently used for
deterministic approximations of permanents \cite{lsw} .(Much more recent proof of {\bf (Bapat-bound)} was actually
motivated by the scaling algorithm \cite{GS} ,\cite{GS1}  to approximate mixed discriminants and mixed volumes.)
It is easy to check that {\bf (Schrijver-bound)} is implied by {\bf (VDW-bound)} for $k \geq n$ : $\frac{n!}{n^{n}} = \prod_{1 \leq k \leq n}(\frac{k-1}{k})^{k-1} > (\frac{k-1}{k})^{(k-1)n}$ .
Therefore , it was not clear whether the scaling algorithm in \cite{lsw} gives better approximating exponent for sparse matrices :
the "scaled" doubly stochastic matrix may have irrational entries even if the input matrix is boolean and {\bf (Schrijver-bound)} is superior
to the {\bf (VDW-bound)} only on "very" rational sparse doubly stochastic matrices .\\
Since  our generalized  {\bf (Schrijver-bounds)} (18),(19)
depend only on the "sparsity" hence the scaling algorithm for permanents in \cite{lsw} indeed gives better approximating exponent for sparse matrices
(scaling algorithm for mixed discriminants in \cite{GS} ,\cite{GS1} gives better approximating exponent for tuples of "small" rank PSD matrices).

\subsection{Van der Waerden / Schrijver-Valiant like conjectures and homogeneous polynomials}
Let $Hom(m,n)$ be the linear space of homogeneous polynomials $p(x), x \in R^m$ of degree $n$ in $m$ real varibles ;
correspondingly $Hom_{+}(m,n) (Hom_{++}(m,n))$ be a subset of homogeneous polynomials $p(x), x \in R^m$ of degree $n$ in $m$ real varibles 
and nonnegative(positive) coefficients .
\dfn
\begin{enumerate}
\item
Let $p \in Hom_{+}(n,n) , p(x_1,...,x_n) =\sum_{ (r_1,...,r_n) \in I_{n,n} } a_{(r_1,...,r_n) } \prod_{1 \leq i \leq n} x_{i}^{r_{i}}$ 
be a homogeneous polynomial of degree $n$ in $n$ real variables. Here $I_{m,n}$ stands for the set of vectors
$r = (r_1,...,r_m)$ with nonnegative integer components and  $\sum_{1 \leq i \leq m} r_i = n$.\\
The support of the polynomial $p(x_1,...,x_n)$ as above 
is defined as $supp(p) = \{(r_1,...,r_n) \in I_{n,n} : a_{(r_1,...,r_n)} \neq 0 \}$ . The convex hull $CO(supp(p))$ of $supp(p)$ is called
the Newton polytope of $p$ .\\
For a subset $A \subset \{1,...,n\}$ we define $S_{p}(A) = \max_{(r_1,...,r_n) \in supp(p)} \sum_{i \in A} r_i$.
Given a vector $(a_1,...,a_n)$ with positive real coordinates , consider univariate polynomials
$D_{A}(t) = p(t (\sum_{i \in A} e_i) + \sum_{1 \leq j \leq n} a_{j}e_j) , V_{A}(t) =p(t (\sum_{i \in A} e_i) + \sum_{j \in A^{\prime}} a_{j}e_j) $ . Then $S_{p}(A)$ is equal to
the degree of the polynomials $D_{A} ,V_{A}(t) $ :
\beqn
S_{p}(A) = deg(D_{A}) = deg(V_{A})
\eeqn
\item
The following linear differential operator maps $Hom(n,n)$ onto $Hom(n-1,n-1)$ :
$$
p_{x_{1}}(x_2,...,x_n) = \frac{\partial}{\partial x_1} p(0,x_2,...,x_n) .
$$
We define $p_{x_{i}} , 2 \leq i \leq n$ in the same way for all polynomials $p \in Hom(n,n)$. Notice that
\beqn
p(x_1,...,x_n) = x_{i} p_{x_{i}}(x_2,...,x_n) + q(x_1,...,x_n) ; q_{x_{i}} = 0 .
\eeqn
The following inequality follows straight from the definition :
\beqn
S_{p_{x_1}}(A) \leq \min(n-1,S_{p}(A)) : A \subset \{2,...,n\} , p \in Hom_{+}(n,n) .
\eeqn

\item
Consider $p \in Hom_{+}(n,n)$  We define the {\bf Capacity} as 
$$
Cap(p) = \inf_{x_i > 0, \prod_{1 \leq i \leq n }x_i = 1} p(x_1,...,x_n) .
$$
It follows that if $p \in Hom_{+}(n,n)$ then
\beqn
Cap(p) \geq \frac{\partial^{n}}{\partial x_{1} \cdots \partial x_{n}} p(0,0,...,0)
\eeqn
Notice that  
$$
\log(Cap(p)) = inf_{\sum_{1 \leq i \leq n} y_i = 0} \log(p(e^{y_1},...,e^{y_n})) ,
$$
and if $p \in Hom_{+}(n,n)$ then the functional $\log(p(e^{y_1},...,e^{y_n}))$ is convex .
\item
Consider a stratified set of homogeneous polynomials : $F = \bigcup_{1 \leq n < \infty } F_n$ ,
where $F_n \in Hom_{+}(n,n)$ . We call such set {\bf VDW-FAMILY} if it satisfies the following properties :
\begin{enumerate}
\item
If a polynomial $p \in F_n , n > 1$ then for all $1 \leq i \leq n$ the polynomials $p_{x_{i}} \in F_{n-1}$.
\item 
\beqn 
Cap(p_{x_i}) \geq  g(S_{p}(\{i\})) Cap(p) :  p \in F_j , 1 \leq i \leq j ; g(k) = (\frac{k-1}{k})^{k-1} , k \geq 1 .
\eeqn
\end{enumerate}
\end{enumerate}
\edfn

\xam
Let $A = \{A(i,j) : 1  \leq i \leq n \}$ be $n \times n$ matrix with nonnegative entries . Assume that $\sum_{1 \leq j \leq n} A(i,j) > 0$ for all $1  \leq i \leq n$.
Define the following homogeneous polynomial $Mul_{A}(t_1,...,t_n) = \prod_{1 \leq i \leq n} \sum_{1 \leq j \leq n} A(i,j)t_j$ . Clearly ,
$Mul_{A} \in Hom_{+}(n,n)$ and $Mul_{A} \neq 0$ . It is easy to check that $S_{Mul_{A}}(\{j\}) = |\{i : A(i,j) \neq 0 \}|$
($S_{Mul_{A}}(\{j\})$ is equal to the number of non-zero entries in the $j$th column of $A$) . \\
{\bf Notice that if $A \in \Lambda(k,n)$ (or  $A \in \Omega(k,n)$) then $S_{Mul_{A}}(\{j\}) \leq k , 1 \leq j \leq n$ .}\\
More generally , consider a $n$-tuple ${ \bf A } = (A_{1}, A_{2},... A_{n})$ , where the complex hermitian $n \times n$ matrices are positive
semidefinite and $\sum_{1 \leq i \leq n} A_{i} \succ 0$ (their sum is positive definite). Then the homogeneous polynomial
$ DET_{{ \bf A }}(t_1,...,t_n) = \det (\sum_{1 \leq i \leq n} t_{i}A_{i}) \in Hom_{+}(n,n)$ and $ DET_{{ \bf A }}\neq 0$ .\\
Similarly to polynomials $Mul_{A}$ , we get that $S_{DET_{{ \bf A }}}(\{j\}) = Rank(A_j) , 1 \leq j \leq n$. \\
The Van Der Waerden conjecture on permanents as well as Bapat's conjecture on mixed discriminants
can be equivalently stated in the
following way (notice the absence of doubly stochasticity ):
\beqn
\frac{n!}{n^{n}} Cap(q) \leq \frac{\partial^n}{\partial x_1...\partial x_n} q(0,...,0) \leq Cap(q)
\eeqn
The van der Waerden conjecture on the permanents corresponds to polynomials $Mul_{A} \in Hom_{+}(n,n) : A \geq 0$ ,
the Bapat's conjecture on mixed discriminants corresponds to $DET_{{ \bf A }} \in Hom_{+}(n,n) : { \bf A } \succeq 0$ .
The connection between inequality (7) and the standard forms of the van der Waerden and Bapat's conjectures is established with the
help of the scaling (\cite{lsw} , \cite{GS} , \cite{GS1}). Notice that the functional
$\log(p(e^{y_1},...,e^{y_1}))$ is convex if $p \in Hom_{+}(n,n)$. Thus the inequality (7) allows a convex relaxation of
the permanent of nonnegative matrices and the mixed discriminant of semidefinite tuples .  This observation was implicit in \cite{lsw} and crucial in \cite{GS} , \cite{GS1} .

\exam

\subsection{ The Main (polynomial) Idea}
The following (meta)theorem describes the main idea of this paper .
\thm
Let $F = \bigcup_{1 \leq n < \infty } F_n$ be a {\bf VDW-FAMILY}  and the homogeneous polynomial $p \in F_n$.
Then the following inequality holds :
\beqn
\prod_{1 \leq i \leq n} g(\min(S_{p}(\{i\})),n+1-i)) Cap(p) \leq \frac{\partial^n}{\partial x_1...\partial x_n} p(0,...,0) \leq Cap(p) .
\eeqn 
\ethm

\prf
Our proof is by natural induction . Notice that the function $g(k) = (\frac{k-1}{k})^{k-1}$ is strictly decreasing on the semiinterval $[1,\infty)$ 
(we define $g(0)=1$). The theorem is obviously true for $n=1$ . Suppose it is true for all $n \leq k < \infty$ and
the polynomial $p \in F_{k+1}$ . Since $F = \bigcup_{1 \leq n < \infty } F_n$ is a {\bf VDW-FAMILY} hence
$ p_{x_{k+1}} \in F_{k} ; Cap(p_{x_{k+1}}) \geq  g(S_{p}(\{k+1\}))Cap(p)$ . It follows from (obvious) inequality (4) 
that $S_{p_{x_{k+1}}}(\{i \}) \leq S_{p}(\{i \}) , 1 \leq i \leq k$ . Therefore , by induction , we get the needed inequality :\\

$\frac{\partial^{k+1}}{\partial x_1...\partial x_{k+1}}p(0,...,0) = \frac{\partial}{\partial x_{k+1}}(\frac{\partial^{k}}{\partial x_1...\partial x_{k}}p_{x_{k+1}}(0,...,0))
\geq$ \\
$\geq \prod_{1 \leq i \leq k} g(\min(S_{p_{x_{k+1}}}(\{i\})),n+1-i)) Cap(p_{x_{k+1}}) \geq \\
\geq \prod_{1 \leq i \leq k} g(\min(S_{p_{x_{k+1}}}(\{i\})),n+1-i)) Cap(p) g(S_{p}(\{k+1\})) \geq\\
\geq \prod_{1 \leq i \leq k+1} g(\min(S_{p}(\{i\})),n+1-i)) Cap(p) .$
\eprf

\cor
\begin{enumerate}
\item
If the homogeneous polynomial $p \in F_n$ then
\beqn
\frac{n!}{n^n} Cap(p) \leq \frac{\partial^n}{\partial x_1...\partial x_n} p(0,...,0) \leq Cap(p) .
\eeqn
\item
If the homogeneous polynomial $p \in F_n$ and $S_{p}(\{i\})) \leq k \leq n , 1 \leq i \leq n$ then
\beqn
(\frac{k -1}{k})^{(k -1)(n-k)} \frac{k!}{k^{k}} Cap(p) \leq \frac{\partial^n}{\partial x_1...\partial x_n} p(0,...,0) \leq Cap(p) .
\eeqn
\end{enumerate}
\ecor

\prf
Both inequalities follow the main inequality (8) and from the next easily proved identity 
\beqn
\frac{n!}{n^n} = \prod_{1 \leq k \leq n} g(k) .
\eeqn
\eprf

What is left now is to present a {\bf VDW-FAMILY} which contains all polynomials $ DET_{{ \bf A }}$ , where the $n$-tuple ${ \bf A } = (A_1,...,A_n)$
consists of positive semidefinite hermitian matrices (and thus contains all polynomials $Mul_{A}$ , where $A$ is $n \times n$ matrix with nonnegative entries).
If such  {\bf VDW-FAMILY} set exists then the Van der Waerden , Bapat , Schrijver-Valiant conjectures would follow (without any extra work , see Example 1.2) from
Theorem 1.3 and Corollary 1.4  .\\
One of such {\bf VDW-FAMILY} , consisting of $POS$-hyperbolic polynomials , is defined in the next section.

\section{Hyperbolic polynomials}

{\it The following concept  of hyperbolic polynomials was originated in the theory of partial differential equations \cite{gar}, \cite{horm} ,\cite{kry} . 
It recently became "popular" in the optimization literature \cite{gul} ,\cite{lewis},\cite{ren}. The paper \cite{ren} gives nice
and concise introduction to the area (with much simplified proofs of the key theorems) .}
\dfn
\begin{enumerate}
\item
A homogeneous polynomial $p: C^m \rightarrow C$ of degree $n$( $p \in Hom(m,n)$) is called hyperbolic in the direction $e \in R^m$ 
(or $e$- hyperbolic) if $p(e) \neq 0$ and for each vector $X \in R^m$ the univariate (in $\lambda$) polynomial 
 $p(X - \lambda e)$   has exactly $n$ real roots counting their multiplicities. 
\item
Denote an ordered vector of roots of $p(x - \lambda e)$ as 
$\lambda(X) = (\lambda_{n}(X) \geq \lambda_{n-1}(X) \geq ... \lambda_{1}(X)) $. 
Call $X \in R^m$ $e$-positive ($e$-nonnegative) if $\lambda_{1}(X) > 0$ ($\lambda_{n}(X) \geq 0$).
We denote the closed set of $e$-nonnegative vectors as $N_{e}(p)$,
and the open set of $e$-positive vectors
as $C_{e}(p)$.
\end{enumerate}
\edfn

\dfn
Let $p: C^m \rightarrow C$ be a homogeneous polynomial of degree $n$ in $m$ variables.
Following \cite{khov} , we define the $p$-mixed form of an $n$-vector tuple ${\bf X} = (X_1,..,X_n) : X_i \in C^m$ as
\beqn
M_{p}({\bf X}) = : M_{p}(X_1,..,X_n) = \frac{\partial^n}{\partial \alpha_1...\partial \alpha_n} p(\sum_{1 \leq i \leq n} \alpha_i X_i)
\eeqn
The following polarization identity is well known 
\beqn
M_{p}(X_1,..,X_n) = 2^{-n} \sum_{b_{i} \in \{-1, 1 \}, 1 \leq i \leq n }  p(\sum_{1 \leq i \leq n } b_i X_i) \prod_{1 \leq i \leq n } b_i
\eeqn
Associate with any vector $r = (r_1,...,r_n) \in I_{n,n}$ an $n$-tuple of $m$-dimensional vectors  ${\bf X}_{r}$ consisting
of  $r_i$ copies of $x_{i}  (1 \leq i \leq n) $.  
It follows from the Taylor's formula  that
\beqn
p(\sum_{1 \leq i \leq n} \alpha_i X_i) = \sum_{r \in I_{n,n}} \prod_{1 \leq i \leq n } \alpha_{i}^{r_{i}} M_{p}({\bf X}_{r})
\frac{1}{\prod_{1 \leq i \leq n } r_{i}!}
\eeqn
\edfn
We collected in the following proposition the properties of hyperbolic polynomials used in this paper .
\pro
{\bf FACT 1 .}
\beqn
p(X) = p(e) \prod_{1 \leq i \leq n }\lambda_{i}(X) .
\eeqn
{\bf FACT 2 .}
If $p$ is $e$-hyperbolic polynomial and $p(e)$ is a real nonzero number then the coefficients of $p$ are real (\cite{horm} , follows from (15) via the standard
interpolation). If $p$ is $e$-hyperbolic polynomial and $p(e)>0$ then $p(X) > 0$ for all $e$-positive vectors $X\in C_{e}(p) \subset R^m$ . \\

{\bf FACT 3 .}
Let $p \in Hom(m,n)$ be $e$-hyperbolic polynomial and $d \in C_{e}(p) \subset R^m$ . Then $p$ is also $d$- hyperbolic and 
$C_{d}(p) = C_{e}(p) , N_{d}(p) = N_{e}(p) $ . (\cite{gar} , \cite{khov} , very simple proof in \cite{ren} .) \\

{\bf FACT 4 .} Let $p \in Hom(m,n)$ be $e$-hyperbolic polynomial . Then the polynomial $p_{e}(X) =: \frac{d}{dt} p(X+te)|_{(t=0)} ; p_{e} \in Hom(m,n-1)$
is also $e$- hyperbolic and $C_{e}(p) \subset C_{e}(p_{e})$   (\cite{khov} , \cite{ren} , Rolle's theorem ).

{\bf FACT 5 .} Let $p \in Hom(m,n)$. Then the $p$-mixed form $M_{p}(X_1,..,X_n)$ is linear in each vector argument $X_i \in C^{m}$.
Let $p\in Hom(m,n)$ be $e$- hyperbolic and $p(e) > 0$ . Then $M_{p}(X_1,..,X_n) > 0$ if the vectors $X_i \in R^{m} , 1 \leq i \leq n$ are $e$-positive 
(\cite{khov} , proved by induction using FACT 4) .
\epro
We use in this paper the following sub-class of hyperbolic polynomials .
\dfn
A polynomial $p \in Hom(m,n)$ is called $POS$-hyperbolic if $p(e) > 0 , e =(1,1,...,1) \in R^m$ ; $p$ is $e$-hyperbolic
and the closed convex cone $N_{e}(p)$ contains the nonnegative orthant $R^{m}_{+}$ . (In other words ,
all the roots of the univariate polynomial equation $p(X -te) =0$ are nonnegative if the coordinates of the vector
$X$ are nonnegative real numbers.) \\
It follows from the identity (14) and {\bf FACT 5} that $POS$-hyperbolic polynomials have nonnegative coefficients .
\edfn

Probably the best known example of a hyperbolic polynomial comes from the hyperbolic geometry :
$p(x_0 ,...,x_k) =x_0^{2} - \sum_{1 \leq i \leq k } x _i^{2}$.  
This polynomial is hyperbolic in the direction $(1,0,0,...,0)$.  Another "popular" hyperbolic polynomial is
$\det(X)$ restricted on a linear real space of hermitian $n \times n$ matrices .
In this case mixed forms are just mixed discriminants , hyperbolic direction is the identity matrix $I$ ,
the corresponding closed hyperbolic cone of $I$-nonnegative vectors coincides with a closed convex cone of positive semidefinite matrices .\\ 
Less known , but very interesting , hyperbolic polynomial is the Moore determinant $M\det(Y)$ 
restricted on a linear real space of hermitian quaternionic $n \times n$ matrices .
(The Moore determinant is a particular case of the generic norms on Jordan Algebras.)
The Moore determinant is , essentially , the Pfaffian (see the corresponding definitions and
the theory in a very readable paper \cite{quat} ) .\\
{\bf This paper benefits from the fact that as multilinear polynomials $Mul_{A} \in Hom_{+}(n,n) : A \geq 0 , Ae > 0$ ,
as well determinantal polynomials $DET_{{ \bf A }} \in Hom_{+}(n,n) : { \bf A } \succeq 0 , \sum_{1 \leq i \leq n} A_i \succ 0$ are $POS$-hyperbolic.}

\subsection{$POS$-Hyperbolic polynomials  form {\bf VDW-FAMILY} }

Let $q \in Hom_{+}(n,n)$ be a $POS$-hyperbolic polynomial. For a vector $X \in C^{n}$ we define the integer number $Rank_{q}(X)$ as the number
of nonzero roots of the equation $q(X-te)=0 , e =(1,1,...,1) = \sum_{1 \leq i \leq n} e_i$. It follows
from the identity (2) that $Rank_{q}(\sum_{i \in A} e_i) = S_{q}(A) , A \subset \{1,2,...,n\}$ .

\thm
\begin{enumerate}
\item
Let $q \in Hom_{+}(n,n)$ be $POS$-hyperbolic polynomial . If $1 \leq Rank_{q}(e_1) = k \leq n$ then
\beqn
Cap(q_{x_1}) \geq g(k) Cap(q) 
\eeqn
\item
Let $q(x_1,x_2,...,x_n)$ be a $POS$-hyperbolic (homogeneous) polynomial of degree $n$ .
Then either the polynomial $q_{x_1} = 0$ or $q_{x_1}$ is $POS$-hyperbolic . 
If $Cap(q) >0$ then $q_{x_1}$ is (nonzero) $POS$-hyperbolic . 
\end{enumerate}
\ethm
\cor
Let $PHP(n) \subset Hom_{+}(n,n)$ be a set of $POS$-hyperbolic polynomials of degree $n$ in $n$ variables ;
define $PHP_{+}(n) = \{p \in PHP(n) : Cap(p) > 0 \}$. Then as $\cup_{n \geq 1} (PHP(n) \cup \{0\})$ as well
$\cup_{n \geq 1} PHP_{+}(n) $ is {\bf VDW-FAMILY} .
\ecor
The second part of Theorem 2.5 is , up to minor modifications ,  well known ({\bf FACT 4}; see , for instance, \cite{khov}, \cite{ren} )
and follows from Rolle's theorem .
The main new results "responsible" for the first part of Theorem 2.5 are the next Lemma 2.7 and its Corollary 2.8.

\lem
\begin{enumerate}
\item
Let $c_1,...,c_n$ be real numbers ; $0 \leq c_i \leq 1 , 1 \leq i \leq n$ and $\sum_{1 \leq i \leq n} c_i = n-1$.\\
Define the following symmetric functions :
$$
S_{n}  =\prod_{1 \leq i \leq n} c_i , S_{n-1} = \sum_{1 \leq i \leq n} \prod_{j \neq i} c_j .
$$
Then the following entropic inequality holds :
$$
S_{n-1} - n S_{n} \geq e^{\sum_{1 \leq i \leq n} c_i \log(c_{i})} .
$$

\item (Mini van der Waerden conjecture)\\
Consider a doubly-stochastic $n \times n$ matrix $A= [a|b|...|b]$ .
I.e. $A$ has $n-1$ columns equal to the column vector $b$ , and one column equal to the column vector $a$ .
Let $a = (a_1,...,a_n)^{T} : a_i \geq 0 , \sum_{1 \leq i \leq n} a_i = 1$ ;
$b = (b_1,...,b_n)^{T} : b_i = \frac{1 - a_i}{n-1} ,1 \leq i \leq n .$
Then the permanent $Per(A) \geq \frac{n!}{n^{n}}$ .
\end{enumerate}
\elem

\cor
Consider an univariate polynomial \\
$R(t) = \sum_{0 \leq i \leq n} d_i t^{i} = \prod_{1 \leq i \leq n} (a_i t + b_i)$ ,
where $a_i , b_i \geq 0$ . If for some positive real number $C$ the inequality 
$R(t) \geq C t$ holds for all $t \geq 0$ then
\beqn
d_1 = \frac{\partial}{\partial t} R(0) \geq C((\frac{n-1}{n})^{n-1})
\eeqn
The inequality (17) is attained only on the polynomials $R(t) = A (t + a)^{n} : A,a > 0 .$
\ecor

\subsection{Newton Inequalities , Alternative Proof of Corollary 2.8 , Volume Polynomials}
Let $R(t) = \sum_{0 \leq i \leq n} d_i t^i$ be an univariate polynomial with real coefficients .
If such polynomial $R$ has all real roots then its coefficients satisfy the
following Newton's inequalities :
$$
NIs :  d_i^{2} \geq d_{i-1} d_{i+1} \frac{{n \choose i}^2}{{n \choose i-1}{n \choose i+1}} : 1 \leq i \leq n-1 .
$$
The following weak Newton's inequalities $WNIs$ follow from $NIs$ if the coefficients are nonnegative:
$$
WNIs : d_i d_{0}^{i-1} \leq \frac{d_{1}}{n}^{i} {n \choose i} : 2 \leq i \leq n .
$$

\lem
Let $R(t) = \sum_{0 \leq i \leq n} d_i t^i$ be an univariate polynomial with real nonnegative coefficients satisfying
weak Newton's inequalities $WNIs$ .If for some positive real number $C$ the inequality 
$R(t) \geq C t$ holds for all $t \geq 0$ then
$$
d_1 \geq C((\frac{n-1}{n})^{n-1}) .
$$
\elem

\prf
If $d_0 =0$ then $d_1 \geq C > C((\frac{n-1}{n})^{n-1})$ . Thus we can assume that $d_0 = 1$ . It follows from
weak Newton's inequalities $WNIs$ that 
$$
d_i \leq (\frac{d_{1}}{n})^{i} {n \choose i} : 2 \leq i \leq n .
$$
Therefore for nonnegative values of $t \geq 0$ we get the inequality
$$
R(t) \leq 1 + (\frac{d_1 t}{n}) {n \choose 1} + (\frac{d_1 t}{n})^{2} {n \choose 2} +... (\frac{d_1 t}{n})^{n} {n \choose n} = (1 + \frac{d_1 t}{n})^{n} .
$$
Which gives the inequality $(1 + \frac{d_1 t}{n})^{n} \geq Ct$ . The inequality $d_1 \geq C((\frac{n-1}{n})^{n-1})$ follows now easily .
\eprf

\rem
The Newton Inequalities are not sufficient for the real rootedness . The classical example is provided by some univariate
volume polynomials $R(t) = Vol(t C_1 +C_2)$ , where $C_1,C_2$ are convex compact sets (see , for instance , \cite{khov}) .
In this case the Newton Inequalities follow from the celebrated Alexandrov-Fenchel Inequalities. \\
Using Lemma 2.9 , the Alexandrov-Fenchel Inequalities and a bit of extra work allows to extend the results of this
paper , i.e. Theorem 1.3 ,to the multivariate volume polynomials $Vol(t_1 C_1 + t_2 C_2 +...+ t_n C_n)$ , where
$C_1,C_2,...,C_n$ are convex compact subsets of $R^n$ . In other words , there exists a {\bf VDW-FAMILY} which contains
all such volume polynomials $Vol(t_1 C_1 + t_2 C_2 +...+ t_n C_n)$ .\\
This extension leads to a randomized poly-time
algorithm to approximate the mixed volume $M(C_1,...,C_n) = \frac{\partial^n}{\partial t_1...\partial t_n} Vol(t_1 C_1 + t_2 C_2 +...+ t_n C_n)$
within exponential factor $e^n$ . The algorithm is pretty much the same as in Theorem 4.7 in this paper , the randomization
is needed to evaluate the oracle , i.e. to evaluate $Vol(t_1 C_1 + t_2 C_2 +...+ t_n C_n)$ .\\
The best current appoximation factor is $n^{O(n)}$ (\cite{Bar1} , \cite{Bar2} (randomized) ; \cite{GS} , \cite{GS1} (deterministic) .
\erem

\dfn
Let $p \in Hom_{+}(n,n) , p(x_1,...,x_n) =\sum_{ (r_1,...,r_n) \in I_{n,n} } a_{(r_1,...,r_n) } \prod_{1 \leq i \leq n} x_{i}^{r_{i}}$ 
be a homogeneous polynomial with nonnegative coefficients of degree $n$ in $n$ real variables. Call such polynomial
{\bf AF-Polynomial} if the following Alexandrov-Fenchel Inequalities hold :
$$
M_{p}(X_1,X_2,X_3,...,X_n)^{2} \geq M_{p}(X_1,X_1,X_3,...,X_n) M_{p}(X_2,X_2,X_3,...,X_n) : X_1,...,X_n \in R^{n}_{+} .
$$
(The $p$-mixed form $M_{p}(X_1,X_2,X_3,...,X_n)$ is defined in Definition 2.2 (formula (12) .) \\
We denote as $AF(n)$ a set of all {\bf AF-Polynomial} of degree $n$ in $n$ real variables  and define \\
$AF_{+}(n) : \{p \in AF(n) : Cap(p) > 0\}$ .
We denote as $Vol(n)$ a set of polynomials $Vol(t_1 C_1 + t_2 C_2 +...+ t_n C_n)$ , where
$C_1,C_2,...,C_n$ are convex compact subsets of $R^n$  .

{\it Notice that $PHP(n) \subset AF(n)$ ( $POS$-Hyperbolic polynomials are {\bf AF-Polynomials}) \cite{khov} ;
$Vol(n) \subset AF(n)$ , this inclusion is just a restatement of the celebrated Alexandrov-Fenchel Inequalities for mixed volumes \cite{Al38}, \cite{Sch93}. }

\edfn

\thm
as $\cup_{n \geq 1} AF(n) $ as well
$\cup_{n \geq 1} AF_{+}(n) $ is {\bf VDW-FAMILY} .
\ethm

\prf
The definition of the {\bf VDW-FAMILY} consists of two properties (see Part 4 of Definition 1.1) .
The property (a) follows from the definition of {\bf AF-Polynomials} :
$$
p_{x_{1}}(x_2,...,x_n) = \frac{\partial}{\partial x_1} p(0,x_2,...,x_n) = ((n-1)!)^{-1} M_p(e_1,Y,Y,...,Y) ;
Y= (0,x_2,...,x_n) .
$$
The property (b) follows from Lemma 2.9 . Indeed , if $X,Y \in R^{n}_{+}$ and $p$ is a {\bf AF-Polynomial} then
the coefficients of the univariate polynomial $p(tX + Y)$ are nonnegative and satisfy
the Newton inequalities $NIs$  .
\eprf

\section{Harvest}
Corollary 2.6 allows to "plug-in" $POS$-hyperbolic polynomials to Theorem 1.3 . The most spectacular application
is the following generalization of {\bf (Schrijver-bound)} .
\thm
Let $A = \{A(i,j) : 1 \leq i,j \leq n\}$ be a matrix with nonnegative entries .
Define $C_j = Card(\{i : A(i,j) \neq 0 \}) , 1 \leq j \leq n$. I.e. $C_j$ is the number of non-zero entries in the $j$th column of $A$ .
Then
\beqn
per(A) \geq \prod_{1 \leq j \leq n} g(\min(C_{j},n+1-j)) Cap(Mul_{A}) .
\eeqn
If $C_j \leq k , 1 \leq j \leq n$ then
\beqn 
per(A) \geq (\frac{k -1}{k})^{(k -1)(n-k)} \frac{k!}{k^{k}} Cap(Mul_{A}) .
\eeqn
(Recall that if $A$ is doubly stochastic then $Cap(Mul_{A}) = 1$ .)
\ethm
\rem
The lower bound (18) can be viewed as a {\bf NONREGULAR} generalization of {\bf (Schrijver-bound)} ; it "interpolates" between
{\bf (VDW-bound)} ($C_j = n$) and the sparse case ($C_j << n$) . The lower bound (19) is actually sharper than
{\bf (Schrijver-bound)} : $\frac{k!}{k^{k}} = \prod_{1 \leq j \leq k} g(j) > g(k)^{K} = (\frac{k -1}{k})^{(k -1)(k)}$.
\erem
\section{Algorithmic Applications}
Suppose that a $POS$-hyperbolic polynomial 
$$
p(x_1,...,x_n) =\sum_{ \sum_{1 \leq i \leq n} r_{i} =n } a_{(r_1,...,r_n) } \prod_{1 \leq i \leq n} x_{i}^{r_{i}}
$$ 
has nonnegative integer coefficients and is given as an oracle . I.e. we don't have a list coefficients ,
but can evaluate $p(x_1,...,x_n)$ on rational inputs . \\
A deterministic polynomial-time oracle algorithm is any algorithm which
evaluates the given polynomial $p(.)$ at a number of rational vectors $q^{(i)}=(q_{1}^{(i)},...,q_{n}^{(i)})$ which is 
polynomial in $n$ and $\log(p(1,1,..,1))$; these rational vectors $q^{(i)}$ are required to 
have bit-wise complexity which is polynomial in $n$ and $\log(p(1,1,..,1))$ ; and the number of additional auxilary arithmetic 
computations is also polynomial in $n$ and $\log(p(1,1,..,1))$ .\\ 
If the number of oracle calls ( evaluations of the given polynomial $p(.)$) , the number of additional auxilary arithmetic computations
and bit-wise complexity of the rational input vectors $q^{(i)}$
are all polynomial in $n$ (no dependence on $\log(p(1,1,..,1))$ ) then such algorithm is called  deterministic strongly polynomial-time oracle algorithm.\\

The following result was proved in \cite{newhyp}.
\thm
\begin{enumerate}
\item
Let $p \in Hom_{+}(n,n)$ be $POS$-hyperbolic polynomial . Then the function $Rank_{p}(\sum_{i \in A} e_i) = S_{p}(A)$ is
submodular , i.e. $S_{p}(A \cup B) \leq S_{p}(A) + S_{p}(B) - S_{p}(A \cap B) : A,B \subset \{1,2,...,n \}$ .
\item
Consider a nonnegative integer vector  $r =(r_1,...,r_n)  , \sum_{1 \leq i \leq n} r_i = n$ . Then
$r \in supp(p)$ iff $r(S) = \sum_{i \in S} r_i \leq S_{p}(S) : S \subset \{1,2,...,n \} $ .
\end{enumerate}
\ethm

\cor
Let $p \in Hom_{+}(n,n)$ be $POS$-hyperbolic polynomial. Associate with this polynomial $p$ the following
bounded convex polytope : 
$$
SUB_{p} = \{ (x_1,...,x_n) : \sum_{i \in S} x_i \leq S_{p}(S) : S \subset \{1,2,...,n \}  ; \sum_{1 \leq i \leq n} x_i = n  ; x_i \geq 0 , 1 \leq i \leq n\}
$$
Then $SUB_{p}$ is equal to the Newton polytope of $p$ , i.e. $SUB_{p} = CO(supp(p)$. 
\ecor
\prf
The inclusion $CO(supp(p)) \subset SUB_{p}$ follows from the definition of the function $S_{p}$ . 
Since the function $S_{p} : 2^{\{1,2,...,n \}} \rightarrow \{0,1,2,...,n\}$ is submodular and integer valued
 hence the extreme points of the polytope $SUB_{p}$ are integer nonnegative vectors \cite{gr:lo:sc} .
Using the second part of Theorem 4.1 , we conclude that all the extreme points of the polytope $SUB_{p}$ belong to $supp(p)$ .
It follows from the Krein-Milman Theorem that $SUB_{p} \subset CO(supp(p))$.
\eprf

\cor
Given $POS$-hyperbolic polynomial $p \in Hom_{+}(n,n)$ as an oracle , there exists strongly polynomial-time oracle algorithm for the membership problem
as for $supp(p)$ as well for the Newton polytope  $CO(supp(p))$.
\ecor
\prf
Let $X=(x_1,...,x_n)$ be a vector with real nonnegative coordinates , $sum_{1 \leq i \leq n} x_i = n$ . Consider
a function $G_{X}(S) =  S_{p}(S) -\sum_{i \in S} x_i , S \subset \{1,2,...,n \}$. Then $G_{X}$ is submodular and 
$X \in CO(supp(p))$ iff $\min_{S \subset \{1,2,...,n \}} G _{X}(S) \geq 0$ ($X \in supp(p)$ iff $\min_{S \subset \{1,2,...,n \}} G _{X}(S) \geq 0$
and $X$ is integer). In view on the recent results on the minimization of submodular functions we only
need to prove that there exists a strogly polynomial-time oracle algorithm to compute $G_{X}(S)$ . Computing $\sum_{i \in S} x_i$ is easy .
And $S_{p}(S) = deg(D_{A})$ , where the univariate polynomial
$D_{A}(t) = p(t (\sum_{i \in A} e_i) + \sum_{1 \leq j \leq n}e_j)$ . Clearly we can compute the degree $deg(D_{A})$ via the standard interpolation ,
which amounts to at most $n+1$ evaluations of $p$ and $O(n \log(n)^{2}$ arithmetic operations .
\eprf

\rem
Consider $p \in Hom_{+}(4,4) , p = x_{1}x_{2}x_{3}x_{4} + x_{2}^{2} x_{4}^{2}$. 
Then $S_{p}(\{1,2,3\}) = 3 , S_{p}(\{1,2\}) = 2 , S_{p}(\{2,3\}) = 2 , S_{p}(\{2\}) = 2$ and therefore the function $S_{p}$ is not submodular .
It is easy to see that $S_{q}$ is submodular for all $q \in Hom_{+}(3,3)$.  On the other hand ,
consider $p \in Hom_{+}(3,3) , p = x_{1}x_{2}^{2} + x_{2}x_{3}^{2} + x_{3}x_{1}^{2}$. One can check that $SUB_{p} \neq CO(supp(p)$.
The proof in \cite{newhyp} of submodularity of $S_{p}$ for $POS$-hyperbolic polynomials $p$ is based on the proved LAX conjecture (\cite{lax} , \cite{hel}, \cite{vin}), very nonelementary result . \\
It was proved in \cite{newhyp} that , unless {\bf P = NP} , there is no deterministic polynomial-time oracle algorithm to check
if $(1,1,...,1) \in supp(p)$ for integer polynomials $p \in Hom_{+}(n,n)$ .  
\erem

\dfn
A homogeneous polynomial $q \in Hom_{+}(n,n)$  is called doubly-stochastic
if  $\frac{\partial}{\partial x_{i } } q(1,1,...,1) = 1$ for all $1 \leq i \leq n$.
The doubly-stochastic defect of the polynomial $q$ is defined as \\
$DS(q) = \sum_{1 \leq i \leq n} (\frac{\partial}{\partial x_{i } } q(1,1,...,1) - 1)^{2}$.\\
A polynomial $q \in Hom_{+}(n,n)$ is called scalable if there exists a positive vector $\beta =(\beta_1,...,\beta_n)$ such that the scaled
polynomial $q_{\beta}(x_1,...,x_n)) = q(\sum_{1 \leq i \leq n} \beta_i x_i)$ is doubly-stochastic .
(It is easy to see that $q \in Hom_{+}(n,n)$ is scalable iff the infimum $\inf_{x_i > 0, \prod_{1 \leq i \leq n }x_i = 1} q(x_1,...,x_n) = Cap(q)$ is
attained .)
A polynomial $q \in Hom_{+}(n,n)$ is called indecomposable if infimum $\inf_{x_i > 0, \prod_{1 \leq i \leq n }x_i = 1} q(x_1,...,x_n) = Cap(q)$
is attained and unique . (Theorem D.1 in Appendix D "justifies" , in the $POS$-hyperbolic case , our notion of indecomposability.)
\edfn

\thm
A $POS$-hyperbolic polynomial $q \in Hom_{+}(n,n)$ is indecomposable if and only if the following two equivalent conditions hold : \\
{\bf Condition 1.}
$\frac{\partial^n}{\partial x_i \partial x_j \prod_{m \neq (i,j)} \partial x_m} q(0,...,0) > 0 : 1 \leq i \neq j \leq n$ . \\
{\bf Condition 2.}
$Rank_{p}(\sum_{i \in A} e_i) = S_{p}(A) > |A| : A \subset \{1,2,...,n \} , 1 \leq |A| < n$ .

\ethm

The following theorem combines the algorithm and its analysis from \cite{GS1} (see section 4 \cite{GS1})in  and Theorem  .
Similarly to  \cite{GS1} , we use the ellipsoid method to approximate $\min_{\sum_{1 \leq i \leq n} \alpha_i = 0} \log(p(e^{\alpha_{1}} ,...,e^{\alpha_{n}})$ .
The starting ball is centered at $0$ and has the radius $n^{\frac{1}{2}} \log(2p(1,1,...,1))$ . To run the ellipsoid method we need to compute
the gradient of $\log(p(e^{\alpha_{1}} ,...,e^{\alpha_{n}})$ ; since $p$ is homogeneous polynomial of degree $n$ hence we can compute the gradient
by $n$ standard univariate interpolations. These univariate interpolations amount to $O(n^2)$ oracle calls and $O(n^2)$ arithmetic operations .
The ellipsoid updating also requires $O(n^2)$ arithmetic operations .

\thm
There exists a deterministic polynomial-time oracle algorithm which computes for given as an oracle indecomposable
$POS$-hyperbolic polynomial $p(x_1,...,x_n)$ a number $F(p)$ satisfying the inequality
$$
\frac{\partial^n}{\partial x_1...\partial x_n} p(0,...,0) \leq F(p) \leq 2(\prod_{1 \leq i \leq n} g(\min(S_{p}(\{i\})),n+1-i)))^{-1}
 \frac{\partial^n}{\partial x_1...\partial x_n} p(0,...,0) \leq 
$$
$$ 
\leq 2\frac{n^n}{n!}\frac{\partial^n}{\partial x_1...\partial x_n} p(0,...,0) .
$$
\ethm

Theorem 4.7 can be (slightly) improved . I.e. it can be applied
to the polynomial 
$$
p_{k}(x_{k+1},...,x_{n}) =\frac{\partial^k}{\partial x_1...\partial x_k} p(0,..,0, x_{k+1},...,x_{n}).
$$
Notice that the polynomial $p_{k}$ is a homogeneous polynomial of degree $n-k$ in $n-k$ variables .
It follows from Theorem 2.5 that if $p= p_{0}$ is $POS$-hyperbolic and $Cap(p) > 0$ then for all $0 \leq k \leq n$
the polynomials $p_{k}$ are also $POS$-hyperbolic and $Cap(p_{k}) > 0$ . 
Also , if $p= p_{0}$ is indecomposable then  $p_{k}$ is indecomposable as well (Theorem 4.6).\\
The trick is that if $k = m \log_{2}(n)$ then (using the polarizational formula (13) )
the polynomial $p_{k}$ can be evaluated using $O(n^{m+1})$ oracle calls of the (original) polynomial
$p= p_{0}$ . This observations allows to decrease the worst case multiplicative factor in Theorem 4.7 from
$e^{n}$ to $\frac{e^{n}}{n^{m}}$ for any fixed $m$ . If the polynomial $p= p_{0}$ can
be explicitly evaluated in deterministic polynomial time , this observation
results in deterministic polynomial time algorithms to approximate
$\frac{\partial^n}{\partial x_1...\partial x_n} p(0,...,0)$ within multiplicative
factor $\frac{e^{n}}{n^{m}}$ for any fixed $m$ . Which is an improvement
of results in \cite{lsw} (permanents , $p$ is a multilinear polynomial) and in \cite{GS} , \cite{GS1} (mixed discriminants,
$p$ is a determinantal polynomial) .

\rem
Let $A$ be an $n \times n$ matrix with nonnegative entries , $S \subset \{1,2,...,n \} , |S| = m\log_{2}(n)$ . Assume , modulo polynomial
time preprocessing , that $A$ is fully indecomposable \cite{lsw} .
Using the Laplace expansion for permanents ,
we get that $per(A) = \sum_{|T| = n -|S|} per(A_{S,T}) per(A_{S^{\prime} , T^{\prime}})$ . This suggest the following deterministic algorithm :
compute exactly the permanents $per(A_{S,T})$ of "small" matrices $A_{S,T}$ and run the algorithm from \cite{lsw} to approximate
with the multiplicative factor $\frac{e^{n}}{n^{m}}$ the permanents of "large" matrices $A_{S^{\prime} , T^{\prime}}$ . This algorithm
achieves the multiplicative factor $\frac{e^{n}}{n^{m}}$ , but it runs in quasi-polynomial time . Our approach
is to apply Theorem 4.7 to indecomposable  $POS$-hyperbolic polynomial $\sum_{|T| = n -|S|} per(A_{S,T}) Mul_{A_{S^{\prime} , T^{\prime}}}$ ,
which can be evaluated in deterministic $Poly(n)$-time. Our new "hyperbolic" {\bf (VDW-bound)} (9) allows multiplicative factor $\frac{e^{n}}{n^{m}}$ .
We can use the same trick for sparse matrices using our new "hyperbolic" {\bf (Schrijver-bounds)} (8),(10).
\erem

\section{Conclusion and Acknowledgements}
Univariate polynomials with nonnegative real roots appear quite often in modern combinatorics , especially in the context
of integer polytopes . The closest to our approach is the class of univariate {\it rook } polynomials \cite{wanless} .
We discovered in this
paper rather unexpected and very likely far-reaching connections between hyperbolic multivariate polynomials and many classical combinatorial and
algorithmic problems .\\
The main "spring" of our approach is that the class
of $POS$-hyperbolic polynomials is large enough to allow the
easy induction . The reader might be surprised by the absence
of Alexandrov-Fenchel inequalities and other ingredients
of proofs in \cite{fal} , \cite{ego} , \cite{gur}.  
In fact , the clearest (in our opinion) proof of the
Alexandrov-Fenchel inequalities for mixed discriminants is
in A.G. Khovanskii' 1984 paper \cite{khov} . The Khovanskii' proof
is based on the similar induction (via partial differentiations) to
the one used in this paper . In a way , the Alexandrov-Fenchel inequalities 
are "hidden" in our proof . \\
Let us summarize the main ideas of our approach :\\
{\bf IDEA 1 .} To facilitate the induction we deal not with doubly stochastic matrices/tuples/polynomials but rather with the {\bf CAPACITY}  
of homogeneous polynomials with nonnegative coefficients .\\
{\bf IDEA 2 .} The notion of the {\bf VDW-FAMILY} allowed to reformulate Van der Waerden / Schrijver-Valiant/Bapat conjectures in
terms of homogeneous polynomials with nonnegative coefficients .\\
{\bf IDEA 3 .} The notion from the theory of linear PDE , $POS$-hyperbolic polynomials , happened to give the needed , i.e. containing
multilinear and determinantal polynomials , {\bf VDW-FAMILY} . Corollary 2.8 , a particularly easy case of
the Van der Waerden Conjecture , was the final strike . \\

I would like to thank Sergey Fomin , Shmuel Friedland , Mihai Putinar, Alex Samorodnitsky , Lex Schrijver , Warren Smith , George Soules , Ian Wanless , Hugo Woerdeman  for
the interest to this paper .\\
This paper was , to a great degree , Internet assisted (sometimes in mysterious ways). Thanks to Google !

\appendix
\section{Proof of Lemma 2.7}
\prf
\begin{enumerate}
\item
Doing simple "algebra" we get that
$$
S_{n-1} - n S_{n} = \prod_{1 \leq i \leq n} c_i (\sum_{1 \leq i \leq n} \frac{1-c_{i}}{c_i} ).
$$
Notice that $0 \leq 1-c_{i} \leq 1$ and $\sum_{1 \leq i \leq n} (1-c_i) = 1$.
Using the concavity of the logarithm we get that
\beqn
\log(S_{n-1} - n S_{n}) \geq \sum_{1 \leq i \leq n} \log(c_i) + \sum_{1 \leq i \leq n} (1-c_{i}) \log(\frac{1}{c_i}) = 
\sum_{1 \leq i \leq n} c_i \log(c_{i}).
\eeqn
\item
$$
per(A) = \frac{ (n-1)!}{ (n-1)^{n-1}} \sum_{1 \leq i \leq n} a_i \prod_{j \neq i} (1-a_j) .
$$
Define $c_i = 1-a_i$ . Then $0 \leq 1-c_{i} \leq 1$ , $\sum_{1 \leq i \leq n} c_i = n-1$ and the permanent
$$
per(A) = \frac{ (n-1)!}{ (n-1)^{n-1}} (S_{n-1} - n S_{n} ) .
$$
It is easy to prove and well known that
$$
\min_{0 \leq 1-c_{i} \leq 1 ;\sum_{1 \leq i \leq n} c_i = n-1 } \sum_{1 \leq i \leq n} c_i \log(c_{i})  = \sum_{1 \leq i \leq n} 
\frac{n-1}{n} \log (\frac{n-1}{n}) =
\log ((\frac{n-1}{n})^{n-1}) .
$$
Using the entropic inequality (20) from the first part we get the following equality
$$
\min_{0 \leq 1-c_{i} \leq 1 ;\sum_{1 \leq i \leq n} c_i = n-1 } S_{n-1} - nS_{n} =
(\frac{n-1}{n})^{n-1} .
$$

Which gives the needed inequality
\beqn
Per(A) \geq \frac{ (n-1)!}{ (n-1)^{n-1}} (\frac{n-1}{n})^{n-1} = \frac{n!}{n^{n}} .
\eeqn

It is easy to see (strict concavity of $\sum_{1 \leq i \leq n} c_i \log(c_{i})$)
that the last inequality is strict unless $ A(i,j) = \frac{1}{n} ; 1 \leq i,j \leq n$.
\end{enumerate}
\eprf
\section{Proof of Corollary 2.8 }
\prf
We can assume WLOG that $a_i > 0 , 1 \leq i \leq n$.If $R(0) = 0$ then
clearly $d_1 = R^{\prime}(0) \geq C \geq C((\frac{n-1}{n})^{n-1})$ . Let $R(0) > 0$ ,
i.e. $b_i > 0 , 1 \leq i \leq n$.\\
Associate with polynomial $R(t) = \prod_{1 \leq i \leq n} (a_i t + b_i)$ the following
matrix with positive entries $A = [a |c|...|c] $ , where $a = (a_1,...,a_n)^{T} , c=\frac{1}{n-1} (b_1,...,b_n)^{T}$.
The condition $R(t) \geq C t , \forall t \geq 0$ is equivalent to the inequality
$Cap(Mul_{A}) \geq C$ . And $per(A) = \frac{(n-1)!} {(n-1)^{n-1}} d_1 .$ 
Since $A$ has positive entries hence there exist two positive diagonal matrices $D_1 , D_2$ such that
$A = D_1 B D_2$ , where the matrix $B$ is doubly stochastic and $B= [f|d|...|d]$ (Sinkhorn Diagonal Scaling ,
see Lemma 3.6 in \cite{GS1} in a more general setting). Since the matrix $B$ is doubly stochastic hence
$Cap(Mul_{B}) =1$ . Thus $Cap(Mul_{A}) = \det(D_1 D_2) Cap(Mul_{B}) = \det(D_1 D_2)$
and $per(A) = \det(D_1 D_2) per(B)$. Therefore , we get from inequality (21) that $per(A) \geq \frac{n!}{n^{n}} Cap(Mul_{A}) \geq \frac{n!}{n^{n}} C$. 

Finally ,it follows that
$$
d_1 = (\frac{(n-1)!} {(n-1)^{n-1}} )^{-1} Per(A) \geq (\frac{(n-1)!} {(n-1)^{n-1}})^{-1} (\frac{n!}{n^{n}} C) = C((\frac{n-1}{n})^{n-1}).
$$
\eprf
\section{Proof of Theorem 2.5  }
We need the following simple result .
\pro
Let $p(X)$ be $e$-hyperbolic (homogeneous) polynomial of degree $n$ , $p(e) > 0$ . 
Consider two $e$-nonnegative vectors
$Z,Y \in N_{e}(p)$ such that $Z+Y \in C_{e}(p)$ , i.e. $Z+Y$ is $e$-positive .
Then
\beqn
p(tZ +Y) = \prod_{1 \leq i \leq n} (a_i t + b_i ) ; a_i,b_i \geq 0 , a_i + b_i > 0 , 1 \leq i \leq n .
\eeqn
\epro

\prf
As the vector $Z+Y = D$ is $e$-positive hence $p(Z+Y) > 0$ ({\bf FACT 1}) ,the polynomial $p$ is $Z+Y$-hyperbolic ({\bf FACT 3}) 
and any $e$-positive ($e$-nonnegative) is also $Z+Y$-positive($Z+Y$-nonnegative) ({\bf FACT 3})  .
Doing simple algebra , we get that $p(tZ +Y) = p((t-1)Z +D)$. \\

Let $0 \leq \lambda_{1} \leq \lambda_{2} \leq ...\leq \lambda_{n} $ be nonnegative roots
of the equation $p(Z- x D) = 0$. Since $D-X = Y \in N_{e}(p) = N_{D}(p) $ hence $\lambda_{n} \leq 1$.
Therefore
$$
p(tZ +Y) = p((t-1)Z +D) = p(D) \prod_{1 \leq i \leq n} (t \lambda_{i} + (1-\lambda_{i}) )
$$
We can put $a_i = (p(Z+Y)) \lambda_{i}  \geq 0 , b_i = (p(Z+Y)) (1-\lambda_{i}) \geq 0$.
\eprf

\subsection{ Proof of the first part of Theorem 2.5 - inequality (16)}
\prf 
Let $q \in Hom_{+}(n,n)$ be $POS$-hyperbolic polynomial and $1 \leq Rank_{q}(e_{1}) = S_{q}(\{1\}) = k$.
Fix positive real numbers $(x_2,...,x_n)$ such that $\prod_{2 \leq i \leq n} x_i = 1.$
Define the following two real $n$-dimensional vectors with nonnegative coordinates :
$Z= (1,0,0,...,0) , Y = (0,x_2,...,x_n)$ . The vector $Z+Y$ is $e$-positive .
Consider the next univariate polynomial $R(t) = q(tZ +Y)$.
It follows from Proposition C.1 that
$$
R(t) = \prod_{1 \leq i \leq n} (a_i t + b_i) = \sum_{0 \leq i \leq n} d_i t^{i} ,
$$
where $a_i , b_i \geq 0$  and $q_{x_{1}}(x_2,...,x_n)= d_1$ and
the cardinality $|\{i : a_i > 0 \}| = k$ (see also equality (2)).  In other words the degree
$deg(R) = k$\\
We get straight from the definition of $Cap(q)$ that 
$$
R(t) = \prod_{1 \leq i \leq n} (a_i t + b_i) = p(t,x_2,...,x_n) \geq Cap(q) t\prod_{2 \leq i \leq n} x_i  = t Cap(q) .
$$
Using Corollary 2.8 , we get that 
$$
q_{x_{1}}(x_2,...,x_n)= d_1 \geq ((\frac{k-1}{k})^{k-1} )Cap(q) .
$$
In other words , that $Cap(r) \geq (\frac{k-1}{k})^{k-1} )Cap(q)$.
\eprf
\subsection{ Proof of the second part }
The second part of Theorem 2.5 is an easy modification of {\bf FACT 4}. We need only to consider
the case $q_{x_{1}} \neq 0$.
We need the following well known fact .
\fac
Consider a sequence of univariate polynomials of the same degree $n$ : $P_{k}(t) = \sum_{0 \leq i \leq n} a_{i,k} t^{i}$ .
suppose that $ \lim_{k \rightarrow \infty} a_{i,k} = a_{i} ,0 \leq i \leq n $ and $a_{n} \neq 0$ .\\
Define $P(t) = \sum_{0 \leq i \leq n} a_{i} t^{i}$ . Then roots of $P_{k}$ converge to roots of $P$ .
In particular if roots of all polynomials  $P_{k}$ are real then also roots of $P$ are real ;
if roots of all polynomials  $P_{k}$ are real nonnegative numbers then also roots of $P$ are real nonnegative numbers .
\efac
It follows from Definition 2.2 and the Taylor's formula (14) that
\beqn
q_{d}(X) =  \frac{d}{dt} q(X+td)|_{(t=0)} = ((n-1)!)^{-1} M_{q}(d,X,X,...,X) : d,X \in R^n 
\eeqn
Notice that $q_{x_{1}}(x_2,...,x_n) = q_{e_{1}}(0,x_2,...,x_n)$ . Consider the following perturbed univariate polynomials $P_{\epsilon}(t) = \sum_{0 \leq i \leq n-1} a_{\epsilon,i} t^{i}$:

$$
P_{\epsilon}(t) = ((n-1)!)^{-1} M_{q}(e_1 + \epsilon e , Y - t((e - e_1) + \epsilon e),...,Y - t((e - e_1) + \epsilon e)) :
$$
$$
e = \sum_{1 \leq i \leq n}e_i , Y = (x_2,...,x_n) \in R^{n-1} , \epsilon > 0 .
$$

We get by a direct inspection that
$$
\lim_{\epsilon \rightarrow \infty} P_{\epsilon}(t) = ((n-1)!)^{-1} M_{q}(e_1 ,Y -t(e - e_1) ,...,Y -t(e - e_1) ) = q_{x_{1}}(x_2-t,...,x_n-t) ; 
$$
$$
\lim_{\epsilon \rightarrow \infty} a_{\epsilon,n-1} = (-1)^{n-1} q_{x_{1}} (1,1,...,1) .
$$
As $q_{x_{1}} \in M_{+}(n-1,n-1)$ and $q_{x_{1}} \neq 0$ hence $q_{x_{1}} (1,1,...,1) > 0$ . Therefore $\lim_{\epsilon \rightarrow \infty} a_{\epsilon,n-1} \neq 0$ .
Since the polynomial $q \in Hom_{+}(n,n)$ is $POS$-hyperbolic hence it follows from {\bf FACT 3} that all the roots of the equation $P_{\epsilon}(t) = 0$ are real ;
if $x_i \geq 0 , 2 \leq i \leq n$ then all the roots are nonnegative . We conclude, using Fact C.2 , that
\begin{enumerate}
\item 
If $Y = (x_2,...,x_n) \in R^{n-1}$ then all the roots of the equation $q_{x_{1}}(Y -t(e - e_1)) = 0$ are real and  $q_{x_{1}}(e - e_1) >0$.
\item 
If $x_i \geq 0 , 2 \leq i \leq n$ then all the roots are nonnegative .
\item
 Since $q_{x_{1}} \in M_{+}(n-1,n-1)$ and $q_{x_{1}} \neq 0$ hence $q_{x_{1}} (x_2,...,x_n) > 0$ if $x_i > 0 , 2 \leq i \leq n$ . It follows from
the equality (15) that all the roots of the equation $q_{x_{1}}(Y -t(e - e_1)) = 0$ are positive if $x_i > 0 , 2 \leq i \leq n$ .
\item 
The polynomial $q_{x_{1}} \in M_{+}(n-1,n-1)$ is $POS$-hyperbolic .
\end{enumerate}

\section{Proof of Theorem 4.6}
\prf
Let $q \in Hom_{+}(n,n)$ be $POS$-hyperbolic polynomial and a pair of indeces $(i \neq j) \subset \{1,2,...,n\}$ .
Define the following integer vectors $r^{(i,j)} = e + e_i - e_j$. {\bf Condition 1} states that $r^{(i,j)} \in supp(q)$ for all such pairs .
The equivalence of {\bf Condition 1} and {\bf Condition 2} follows from the second part of Theorem 4.1.\\
The fact that {\bf Condition 1} implies indecomposability is valid for all polynomials in $ Hom_{+}(n,n)$ and is proved in \cite{GS1}.\\
Suppose that there exists a positive vector $\beta =(\beta_1,...,\beta_n) , \prod_{1 \leq i \leq n}\beta_i =1 $ such that 
$$
q(\beta_1,...,\beta_n) = \inf_{x_i > 0, \prod_{1 \leq i \leq n }x_i = 1} q(x_1,...,x_n) = Cap(q) .
$$
Then the 
polynomial $Q(x_1,...,x_n)) = \frac{1}{q(\beta_1,...,\beta_n)} q(\sum_{1 \leq i \leq n} \beta_i x_i)$ is doubly-stochastic .
Notice that the  $\inf_{x_i > 0, \prod_{1 \leq i \leq n }x_i = 1} q(x_1,...,x_n)$ is attained and unique if and only if
the  $\inf_{x_i > 0, \prod_{1 \leq i \leq n }x_i = 1} Q(x_1,...,x_n)$ is attained and unique ;
it follows from the Euler's identity that $\inf_{x_i > 0, \prod_{1 \leq i \leq n }x_i = 1} Q(x_1,...,x_n) = 1 = Q(e) $ .\\
Assume that {\bf Condition 2} does not hold : there exists a subset $A \subset \{1,2,...,n\} , 1 \leq |A| = m < n$ such that 
$Rank_{q}(\sum_{i \in A} e_i) = Rank_{Q}(\sum_{i \in A} e_i) = m$. Define $e_{A} = \sum_{i \in A} e_i , e_{A^{\prime}} = \sum_{j \in A^{\prime}} e_j = e - e_{A}$ .
Let $\lambda_n \geq ... \geq \lambda_{n-m+1} \geq 0 \geq ...  \geq 0$ be the ordered roots of the equation $Q(e_{A} - te) = 0$ .
Since $e_{A} = e - e_{A^{\prime}}$ and the vector $e_{A^{\prime}}$ is $e$-nonnegative hence $ 0 \leq \lambda_{i} \leq 1 , 1 \leq i \leq n$.
As the polynomial $Q \in Hom_{+}(n,n)$ is doubly-stochastic hence (see \cite{hyp} , \cite{newhyp})  
$$
\sum_{1 \leq i \leq n} \lambda_{i} = \sum_{n-m+1 \leq i \leq n} \lambda_{i} = m = |A| .
$$
Therefore $\lambda_{i} = 1 , n-m+1 \leq i \leq n$ and $\lambda_{j} = 0 , 1 \leq j \leq n-m$ . It follows from the identity (15) that
$Q(a e_{A} + b e_{A^{\prime}}) = 1 = Q(e)$ iff $a^{m}b^{n-m} =1$ , which proves the non-uniqueness of
$\inf_{x_i > 0, \prod_{1 \leq i \leq n }x_i = 1} Q(x_1,...,x_n)$ .
\eprf

The following result is proved very similarly to the previous proof , the only new ingredient is the subadditivity of $Rank_{q}(X) , X \in R^n$.
It "justifies" the notion of "indecomposibility of $POS$-hyperbolic polynomials" . 

\thm
Let $q \in Hom_{+}(n,n)$ be $POS$-hyperbolic polynomial . Supposed that \\
$\inf_{x_i > 0, \prod_{1 \leq i \leq n }x_i = 1} q(x_1,...,x_n)$ is attained
and $Rank_{q}(\sum_{1 \leq i \leq m} e_i) = m , m < n$. \\
Then the polynomial $q$ can be decomposed in the following way :
\beqn
q(x_1,..,x_m,...,x_n) = q_{1}(x_1,..,x_m) q_{2}(x_{m+1},...,x_{n}) :  q_1 \in Hom_{+}(m,m) , q_2 \in Hom_{+}(n-m,n-m)
\eeqn  
\ethm
\end{document}